\documentclass[10pt,twocolumn]{article}          
\usepackage{ijs}                              
\usepackage{times}                               
\usepackage{amsmath,amssymb,latexsym}
\usepackage{fancyhdr,graphicx}
\usepackage[square]{natbib}

\newcommand{\mn}{\mathbb{N}} 
\newcommand{\mz}{\mathbb{Z}} 
\newcommand{\var}[1]{\text{var}\left(#1\right)}
\newcommand{\Esym}{\text{E}}
\newcommand{\E}[1]{\Esym\left[#1\right]}
\newcommand{\Probsym}{\mathbb{P}}
\newcommand{\Prob}[1]{\Probsym\left[#1\right]}

\newtheorem{defn}{Definition}

\begin{document}

\bibliographystyle{plainnat}

\title{A PRACTICAL GUIDE TO MEASURING THE HURST PARAMETER} 
\author{RICHARD G. CLEGG \\ \\
\em Dept. Of Mathematics, University of York, YO10 5DD \\
richard@richardclegg.org }

\date{ \em
\\
This paper describes, in detail, techniques for measuring the Hurst parameter.
Measurements are given on artificial
data both in a raw form and corrupted in various ways to check the robustness
of the tools in question.  Measurements
are also given on real data, both new data sets and well-studied data sets.  
All data and tools used are freely available for download along with simple
``recipes'' which any researcher can follow to replicate these measurements.
}
\affiliation {long-range dependence, data analysis, Hurst parameter, Internet
traffic}

\maketitle
\thispagestyle{fancy} 

\section{INTRODUCTION AND BACKGROUND}

Long-Range Dependence (LRD) is a statistical phenomenon which has received
much attention in the field of telecommunications in the last ten years.
A time-series is described as possessing LRD if it has correlations which
persist over all time scales.  A good guide to LRD is given by
\citep*{beran1994} and a summary in the context of telecommunications is
given by \citep*[chapter one]{clegg2004} (from which some of the material in
this paper is taken).  In the early nineties, LRD was measured in
time-series derived from Internet traffic \citep*{leland1993}.  The importance
of this is that LRD can impact heavily on queuing.  LRD is
characterised by the parameter $H$, the Hurst parameter, (named for a 
hydrologist who pioneered the field in the fifties \citep*{hurst1951}) where
$H \in (1/2,1)$ indicates the presence of LRD.  There are a number of
different statistics which can be used to estimate the Hurst parameter and 
several papers have been written comparing these estimators both in theory
and practice \citep*{taqqu1995,taqqu1997,bardet2003}.
The aim of this paper is not to make a rigorous comparison
of the estimators but, instead, to present a simple and readable guide
to what a researcher can expect from attempting to assess whether LRD
is absent or present in a data set.  All the tools used are available
online using free software.  Software can be downloaded from: \\
{\tt www.richardclegg.org/ \\ lrdsources/software/}

\subsection{Long-Range Dependence in Telecommunications}

In their classic paper, 
Leland et al \citep*{leland1993} measure traffic past a point
on an Ethernet Local Area Network.  They conclude that 
``In the case of Ethernet LAN traffic,
self-similarity is manifested in the absence of a natural length of
a `burst'; at every time scale ranging from a few milliseconds to minutes
and hours, bursts consist of bursty sub-periods separated by less burst
sub-periods.  We also show that the degree of self-similarity (defined via
the Hurst parameter) typically depends on the utilisation level of
the Ethernet and can be used to measure `burstiness' of LAN traffic.''
Since then, a number of authors have replicated these experiments on a
variety of measurments of Internet traffic and the majority
found evidence of LRD or related multi-fractal behaviour.  
Summaries are given in \citep*{sahinoglu1999, willinger2003}.
The reason for the interest in the area is that LRD can, in some
circumstances, negatively impact network performance.  The exact details
of the scale and nature of the effect are uncertain and depend on the
particular LRD process being considered.

\subsection{A Brief Introduction to Long-Range Dependence}

Let $\{X_t: t \in \mn\}$ be a time-series which is weakly stationary (that
is it has a finite mean and the covariance depends only on the separation
or ``lag'' between two points in the series).  Let $\rho(k)$ be the
auto-correlation function (ACF) of $X_t$.  
\begin{defn}
The ACF, $\rho(k)$ for a weakly-stationary time series,  $\{X_t: t \in \mn\}$
is given by
\begin{equation*}
\rho(k)= \frac{\E{(X_t - \mu) (X_{t+k} - \mu)}} {\sigma^2},
\end{equation*}
where $\E{X_t}$ is the expectation of $X_t$, $\mu$ is the mean and $\sigma^2$ is
the variance.
\end{defn}
There are a number of different
definitions of LRD in use in the literature.  A commonly used
definition is given below.
\begin{defn}
The time-series $X_t$ is said to be {\em long-range dependent} if 
$\sum_{k=-\infty}^{\infty} \rho(k)$ diverges.
\label{defn:lrd_weak}
\end{defn}
Often the specific functional form
\begin{equation}
\rho(k) \sim C_\rho k^{-\alpha},
\label{eqn:lrd}
\end{equation}
is assumed where $C_\rho > 0$ and $\alpha \in (0,1)$.  
Note that the symbol $\sim$ is used here and
throughout this paper to mean {\em asymptotically equal to} or 
$f(x) \sim g(x) \Rightarrow f(x)/g(x) = 1$ as $x \rightarrow \infty$ or,
where indicated, as $x \rightarrow 0$.  The parameter $\alpha$ is related
to the Hurst parameter via the equation $\alpha = 2 - 2H$.

If \eqref{eqn:lrd} holds then a similar definition can be shown to hold in the
frequency domain.
\begin{defn}
The {\em spectral density} $f(\lambda)$ of a function with 
ACF $\rho(k)$ and variance $\sigma^2$ can be defined as
\begin{equation*}
f(\lambda)= \frac{\sigma^2}{2 \pi} \sum_{k= -\infty}^{\infty} 
\rho(k) e^{ik\lambda},
\end{equation*}
where $\lambda$ is the frequency, $\sigma^2$ is the variance and $i = \sqrt{-1}$.
\label{defn:spectral}
\end{defn}
Note that this definition of spectral density comes from the
Wiener-Kninchine theorem \citep*{wiener1930}.
\begin{defn}
The weakly-stationary
time-series $X_t$ is said to be {\em long-range dependent} if its
spectral density obeys
\begin{equation*}
f(\lambda) \sim C_f |\lambda|^{-\beta},
\end{equation*}
as $\lambda \rightarrow 0$,
for some $C_f > 0$ and some real $\beta \in (0,1)$.
\label{defn:lrd_freq}
\end{defn}
The parameter $\beta$ is related to the Hurst parameter by 
$H = (1 + \beta)/2$.

LRD relates to a number of other areas of statistics, notably the
presence of statistical self-similarity.  Self-similarity can
be characterised by a self-similarity parameter $H$.
If a self-similar process has stationary increments and $H \in (1/2,1)$
then its increments themselves, taken as a process, form
an LRD process with Hurst parameter $H$.  Indeed analysis
of telecommunications traffic is often described in terms of
self-similarity  and not 
long-range dependence.  (Sometimes the phrase ``asymptotic second-order self-similarity" is used.  This refers to self-similarity in the data when it is
aggregated and is synonymous with LRD.)

In summary, LRD can be thought of in two ways.  In the time domain
it manifests as a high degree of correlation between distantly
separated data points.  In the frequency domain it manifests as a
significant level of power at frequencies near zero.
LRD is, in many ways, a difficult statistical property to work with.
In the time-domain it is measured only at high lags (strictly at
infinite lags) of the ACF --- those very lags where only a few
samples are available and where the measurement errors are largest.
In the frequency domain it is measured at frequencies near zero, again
where it is hardest to make measurements.  Time
series with LRD converge slowly to their mean.  While the Hurst parameter
is perfectly well-defined mathematically, it will be shown that it is,
in fact, a very difficult property to measure in real life.

\section{MEASURING THE HURST PARAMETER}

While the Hurst parameter is perfectly well-defined mathematically,
measuring it is problematic.  The data
must be measured at high lags/low frequencies where fewer readings are
available.  Early estimators were biased and converged
only slowly as the amount of available data increased.  All estimators
are vulnerable to trends in the data, periodicity in the data and
other sources of corruption.  Many estimators assume specific functional
forms for the underlying model and perform poorly if this is misspecified.
The techniques in this paper are chosen for a variety of reasons.  The
R/S parameter, aggregated variance and periodogram are well-known techniques
which have been used for some time in measurements of the Hurst parameter.
The local Whittle and wavelet techniques are newer techniques which
generally fare well in comparative studies.  All the techniques chosen
have freely available code which can be used with free software to
estimate the Hurst parameter.
 
The problems with real-life data are worse than those faced when 
measuring artificial data.  Real life data is likely
to have periodicity (due to, for example, daily usage patterns), trends
and perhaps quantisation effects if readings are taken to a given precision.
The naive researcher taking a data set
and running it through an off-the-shelf method for estimating the
Hurst parameter is likely to
end up with a misleading answer or possibly several different misleading
answers.

\subsection{Data sets to be studied}

A large number of methods are used for generating data exhibiting LRD.
A review of some of the better known methods are given in 
\citep*{bardet2003b}.  In this paper trial data sets with LRD and
a known Hurst parameter are generated using fractional auto-regressive
integrated moving average (FARIMA) modelling and fractional Gaussian
noise (FGN).  The software used to generate the data is included
with at the web address previously mentioned.

A FARIMA model is a well-known time series modelling technique.  It
is a modification of the standard time series
ARIMA ($p,d,q$) model.  
An ARIMA model is defined by
\begin{equation*}
(1 - \sum_{j=1}^p \phi_j \mathbf{B}^j)(1-\mathbf{B})^dX_i=
(1 - \sum_{j=1}^q \theta_j \mathbf{B}^j) \varepsilon_i,
\end{equation*}
where $p$ is the order of the AR part of the model, the $\phi_i$
are the AR parameters, $p$ is the order of the MA part of the model,
the $\theta_j$ are the MA parameters, $d \in \mz$ is the
order of differencing, the $\varepsilon_i$ are i.i.d. noise
(usually normally distributed with zero mean)
and $\mathbf{B}$ is the backshift operator
defined by $\mathbf{B}(X_t) = X_{t-1}$.  If, instead of being
an integer, the model is changed so that $d \in (0,1/2)$ 
then the model is a FARIMA model.  
If the $\phi_i$ and $\theta_i$ are 
chosen so that the model is stationary and $d \in(0,1/2)$
then the model will be LRD with $H= d+ 1/2$.
FARIMA processes were proposed by 
\citep*{granger1980} and a description in
the context of LRD can be found in \citep*[pages 59--66]{beran1994}.

Fractional Brownian Motion is a process $B_H(t)$ for $t \geq 0$ obeying,
\begin{itemize}
\item $B_H(0) = 0$ almost surely,
\item $B_H(t)$ is a continuous function of $t$,
\item The distribution of $B_H(t)$ obeys
\begin{align*}
\Prob{B_H(t+k) - B_H(t) \leq x} =  \\ (2 \pi) ^{-\frac{1}{2} }k ^{-H}
\int \limits_{-\infty}^{x} \exp\left(\frac{-u^2}{2k^{2H}}\right)du,
\end{align*}
\end{itemize}
where $H \in (1/2,1)$ is the Hurst parameter.  The process $B_H(t)$ is
known as fractional Brownian motion (FBM) and its increments are
known as fractional Gaussian noise (FGN).  FBM is a self-similar
process with self-similarity parameter $H$ and, when $H \in (1/2,1)$,
FGN exhibits long-range dependence with Hurst parameter $H$.  
When $H = 1/2$ in
the above, then the process is the well known Weiner process 
(Brownian motion) and the increments are independent (Gaussian noise).
A number of authors
have described computationally efficient methods for generating FGN and
FBM.  The one used in this paper is due to \citep*{paxson1997}.

Data generated from these models will be tested using the various 
measurement techniques and then the same data set will be corrupted in
several ways to see how this disrupts measurements:
\begin{itemize}
\item Addition of zero mean AR(1) model with a high degree of short-range
correlation ($X_t= 0.9 X_{t-1} + \varepsilon_t$).  This simulates a process
with very high local correlations which might be mistaken for a 
long-range dependence.
\item Addition of periodic function (sine wave) --- ten complete 
cycles of a sin wave are added to the signal.  This simulates a 
seasonal effect in the data, for example, a daily usage pattern.
\item Addition of linear trend.  This simulates growth in the data, 
for example the data might be a sample of network traffic at a time 
of day when the network is
growing busier as time continues.
\end{itemize} 
The noise signals are normalised so the standard deviation of 
the corrupting signal
is identical to the standard deviation of the original LRD signal 
to which it is being added.  Note that strictly speaking, while
the addition of an AR(1) model does not change the LRD in the model
and theoretically will leave the Hurst parameter unchanged,
techincally the addition of a trend or of periodic noise makes the
time-series non-stationary and hence the time-series produces are,
strictly speaking, not really LRD.  

In addition, some real-life traffic traces are studied to provide insight
into how well different measurements agree across data sets with and without
various transforms being applied to clean the data.  The data sets used
are listed below.
\begin{itemize}
\item The famous (and much-studied) Bellcore data \citep*{leland1991} 
which was collected in 1989 and has been used for a large number of 
studies since.  Note that, unfortunately, the exact traces used 
in \citep*{leland1993} are
not available for download.  Data from the same sites collected at a similar
time is available online at: \\
{\tt ita.ee.lbl.gov/ \\ html/contrib/BC.html}
\item A data set collected at the University of York in 2001 which consists
of a tcpdump trace of 67 minutes of incoming and outgoing data from
the external link to the university from the rest of the Internet.
\end{itemize}

Three techniques (listed below) were tried to filter real-life 
traces in addition to 
making measurements purely on the raw data.  These methods have been 
selected from the literature as techniques
commonly used by researchers in the
field.  Often a high pass filter would be used to remove
periodicity and trends.  However, since LRD measurements are most
important at low-frequency, that is an obviously inappropriate technique.

\begin{itemize}
\item Transform to log of original data (only appropriate if data is positive).
\item Removal of mean and linear trend (that is, subtract the best fit line
$Y= at  + b$ for constant $a$ and $b$).
\item Removal of high order best-fit polynomial of degree ten (the degree ten was
chosen after higher degrees showed evidence of overfitting).
\end{itemize}

Note that the ``transform to log'' option is not available if the data contains zeros.
In practice some rule of thumb could be considered for replacing zeros with a
minimal value but this substitution was not done here and this pre-processing
technique has not been used where the data contains zeros.

\subsection{Measurement techniques}

The measurement techniques used in this paper can only be described
briefly but references to fuller descriptions with mathematical details are
given.  The techniques used here are chosen for various reasons.  
The R/S statistic, aggregated variance and periodogram are well-known 
techniques with a considerable history of use in estimating long-range
dependence.  The wavelet analysis technique and local Whittle estimator
are newer techniques which perform well in comparative studies and have
strong theoretical backing.

The R/S statistic is a well-known technique for estimating
the Hurst parameter.  It is discussed in \citep*{mandelbrot1969} and 
also \citep*[pages 83--87]{beran1994}.  Let $R(n)$ be the range
of the data aggregated (by simple summation) over blocks of length 
$n$ and $S^2(n)$ be the sample variance of the data aggregated at
the same scale.   
For FGN or FARIMA series the ratio $R/S(n)$ follows
\begin{equation*}
\E{R/S(n)} \sim C_Hn^H,
\end{equation*}
where $C_H$ is a positive, finite constant independent of $n$.
Hence a log-log plot of $R/S(n)$ versus $n$ should have a constant
slope as $n$ becomes large.  A problem with this technique which
is common to many Hurst parameter estimators is knowing which values
of $n$ to consider.  For small $n$ short term correlations dominate
and the readings are not valid.  For large $n$ then there are few
samples and the value of $R/S(n)$ will not be accurate.  Similar 
problems occur for most of the estimators described here.

The aggregated variance technique is described in \citep*[page 92]{beran1994}.
It considers $\var{X^{(m)}}$ where $X^{(m)}_t$ is a time series derrived
from $X_t$ by aggregating it over blocks of size $m$.  The sample variance 
$\var{X^{(m)}}$ should be asymptotically proportional to $m^{2H-2}$ for
large $N/m$ and $m$.

The periodogram, described by \citep*{geweke1983} is defined by
\begin{equation*}
I(\lambda) = \frac{1}{2 \pi N} \left|  
\sum_{j=1}^N X_je^{ij\lambda}\right|^2,
\end{equation*}
where $\lambda$ is the frequency.  For a series with finite variance,
$I(\lambda)$ is an estimate of the spectral density of the series.
From Definition \ref{defn:lrd_freq} then, a log-log plot of $I(\lambda)$
should have a slope of $1 - 2H$ close
to the origin.

Whittle's estimator is a Maxmimum Likelihood Estimator which assumes
a functional form for $I(\lambda)$ and seeks to minimise parameters
based upon this assumption.  
A slight issue with the
Whittle estimator is that the user must specify the functional form
expected, typically either FGN or FARIMA (with the order specified).  If
the user misspecifies the underlying model then errors may occur.
Local Whittle is a semi-parametric version
of this which only assumes a functional form for the spectral density
at frequencies near zero \citep*{robinson1995}. 

Wavelet analysis has been used with success both to measure the
Hurst parameter and also to simulate data \citep*{riedi2003}.  Wavelets
can be thought of as akin to Fourier series but using waveforms other than
sine waves.  The estimator used here fits a straight line to a frequency
spectrum derived using wavelets.  A 95\% confidence interval is given,
however, this should be interpreted only as a confidence interval on
the fitted line and, as will be seen, not as a confidence inteval on
the fitted Hurst parameter.  This is an important distinction --- it is
tempting to consider the confidence intervals given by some estimators
as literally confidence intervals on the measurement of $H$.  Often this
is not the case (as in the case of Wavelet analysis) or is only the case 
if certain conditions are met.

\section{RESULTS}

Results here are in two sections.  Firstly, results are given for simulated
data.  In these cases the expected ``correct'' answer is known and therefore
it can be seen how well the estimators have performed.  The data is then
corrupted by the addition of noise with the same standard deviation as the
original data sets.  Three types of noise are considered as described previously.

In the second section results are given for real data.  The York data is
analysed as a time series of bytes per unit time for two different time units.
The Bellcore data is analysed both in terms of interarrival times and in
terms of bytes per unit time.  Note that, strictly speaking, the interarrival
times do not consititute a proper ``time-series'' since the time units between
readings are not constant.

\subsection{Results on Simulated Data}

For each of the simulation methods chosen, traces have been generated.  Each
trace is 100,000 points of data.  Hurst parameters of 0.7 and 0.9 have been
chosen to represent a low and a high level of long-range dependence in data.
The errors on the wavelet estimator are a 95\% confidence
interval on the fitted regression line (not, as might be thought, the
Hurst parameter measured).

\begin{table*}
\begin{center}
\begin{tabular}{|l| l l l l l|} \hline
Added & R/S Plot & Aggreg. & Period. & Wavelet & Local \\
Noise & & Variance & ogram & Estimate & Whittle \\ \hline
\multicolumn{6}{|c|}{100,000 points FGN --- H= 0.7 --- run one.} \\ \hline
None &0.66  &  0.668  &  0.686  &  0.707 $\pm$ 0.013  &  0.72  \\
AR(1) &0.767  &  0.657  &  0.794  &  0.888 $\pm$ 0.034  &  0.904  \\
Sin &0.667  &  0.969  &  0.692  &  0.707 $\pm$ 0.013  &  0.787  \\
Trend &0.66  &  0.968  &  0.777  &  0.707 $\pm$ 0.013  &  0.766  \\ \hline
\multicolumn{6}{|c|}{100,000 points FGN --- H= 0.7 --- run two.} \\ \hline
None &0.641  &  0.692  &  0.7  &  0.694 $\pm$ 0.007  &  0.721  \\
AR(1) &0.775  &  0.671  &  0.795  &  0.882 $\pm$ 0.036  &  0.902  \\
Sin &0.66  &  0.97  &  0.705  &  0.694 $\pm$ 0.007  &  0.788  \\
Trend &0.641  &  0.968  &  0.769  &  0.694 $\pm$ 0.007  &  0.765  \\ \hline
\multicolumn{6}{|c|}{100,000 points FGN --- H= 0.7 --- run three.} \\ \hline
None &0.636  &  0.69  &  0.704  &  0.708 $\pm$ 0.009  &  0.723  \\
AR(1) &0.734  &  0.654  &  0.79  &  0.876 $\pm$ 0.038  &  0.905  \\
Sin &0.64  &  0.969  &  0.709  &  0.708 $\pm$ 0.009  &  0.787  \\
Trend &0.636  &  0.971  &  0.783  &  0.708 $\pm$ 0.009  &  0.77  \\ \hline
\multicolumn{6}{|c|}{100,000 points FGN --- H= 0.9 --- run one.} \\ \hline
None &0.782  &  0.864  &  0.905  &  0.901 $\pm$ 0.009  &  0.934  \\
AR(1) &0.805  &  0.784  &  0.88  &  0.969 $\pm$ 0.042  &  1.066  \\
Sin &0.772  &  0.961  &  0.907  &  0.901 $\pm$ 0.009  &  0.945  \\
Trend &0.782  &  0.958  &  0.928  &  0.901 $\pm$ 0.009  &  0.939  \\ \hline
\multicolumn{6}{|c|}{100,000 points FGN --- H= 0.9 --- run two.} \\ \hline
None &0.862  &  0.837  &  0.891  &  0.902 $\pm$ 0.003  &  0.933  \\
AR(1) &0.856  &  0.76  &  0.877  &  0.969 $\pm$ 0.038  &  1.062  \\
Sin &0.858  &  0.955  &  0.894  &  0.902 $\pm$ 0.003  &  0.943  \\
Trend &0.862  &  0.954  &  0.921  &  0.902 $\pm$ 0.003  &  0.938  \\ \hline
\multicolumn{6}{|c|}{100,000 points FGN --- H= 0.9 --- run two.} \\ \hline
None &0.793  &  0.884  &  0.907  &  0.904 $\pm$ 0.007  &  0.93  \\
AR(1) &0.818  &  0.802  &  0.871  &  0.972 $\pm$ 0.041  &  1.066  \\
Sin &0.8  &  0.967  &  0.91  &  0.904 $\pm$ 0.007  &  0.943  \\
Trend &0.794  &  0.959  &  0.924  &  0.904 $\pm$ 0.007  &  0.936  \\ \hline
\end{tabular}
\caption{Results for Fractional Gaussian Noise models plus various forms of noise.}
\label{tab:fgn}
\end{center}
\end{table*}

Table \ref{tab:fgn} shows results for various FGN models.  Three runs
each are done with a Hurst parameter of 0.7 and then 0.9.  Firstly it
should be noted that, in all cases, for H=0.7 all estimators are relatively
close when no noise is applied.  The R/S method performs worst, as
it consistently underestimates the Hurst parameter.  The addition of
AR(1) noise confuses all the methods with the Local Whittle performing
particularly poorly.  The correct answer is well outside the confidence
intervals of the Wavelet estimate after this addition (although, as previously
stated, the confidence interval should not be taken literally).  Addition of a
sine wave or a trend causes trouble for the aggregated variance method
but the frequency domain methods (wavelets and local Whittle) do not
seem greatly affected.  

When considering runs with Hurst parameter H=0.9, the R/S method gets
a considerable underestimate even with no corrupting noise.  Note also
that the R/S and aggregated variance method actually produce quite
different estimates for the three runs.  Most
methods seem to perform badly with the AR(1) noise corruption.  Again
the frequency domain methods seem to be able to cope with the sine wave
and with the addition of a trend.  The aggregated variance method seems
to perform particularly badly in the presence of a corrupting sin wave
and a corrupting trend (perhaps not surprising as such a series is no
longer weakly stationary).

\begin{table*}
\begin{center}
\begin{tabular}{|l| l l l l l|} \hline
Added & R/S Plot & Aggreg. & Period. & Wavelet & Local \\
Noise & & Variance & ogram & Estimate & Whittle \\ \hline
\multicolumn{6}
{|c|}{100,000 points FARIMA (0,d,0) --- H = 0.7 --- run one.} \\ \hline 
None & 0.663  &  0.692  &  0.699  &  0.696 $\pm$ 0.004  &  0.681  \\
AR(1) & 0.823  &  0.673  &  0.792  &  0.896 $\pm$ 0.033  &  0.876  \\
Sin & 0.665  &  0.972  &  0.704  &  0.696 $\pm$ 0.004  &  0.765  \\
Trend & 0.662  &  0.973  &  0.786  &  0.696 $\pm$ 0.004  &  0.746  \\ \hline
\multicolumn{6}
{|c|}{100,000 points FARIMA (0,d,0) --- H= 0.7 --- run two.} \\ \hline
None & 0.706  &  0.701  &  0.71  &  0.702 $\pm$ 0.007  &  0.679  \\
AR(1) & 0.837  &  0.673  &  0.791  &  0.891 $\pm$ 0.034  &  0.873  \\
Sin & 0.714  &  0.972  &  0.714  &  0.702 $\pm$ 0.007  &  0.764  \\
Trend & 0.706  &  0.972  &  0.782  &  0.702 $\pm$ 0.007  &  0.742  \\ \hline
\multicolumn{6}{|c|}
{100,000 points FARIMA (0,d,0) --- H= 0.7 --- run three.} \\ \hline
None & 0.718  &  0.684  &  0.696  &  0.687 $\pm$ 0.005  &  0.679  \\
AR(1) & 0.827  &  0.667  &  0.776  &  0.868 $\pm$ 0.044  &  0.872  \\
Sin & 0.723  &  0.973  &  0.701  &  0.687 $\pm$ 0.005  &  0.765  \\
Trend & 0.718  &  0.972  &  0.778  &  0.687 $\pm$ 0.005  &  0.743  \\ \hline
\multicolumn{6}{|c|}{100,000 points FARIMA (1,d,1) --- H= 0.7,
$\phi_1 = 0.5, \theta_1 = 0.5$.} \\ \hline
None & 0.684  &  0.693  &  0.706  &  0.697 $\pm$ 0.006  &  0.68  \\
AR(1) & 0.818  &  0.656  &  0.774  &  0.88 $\pm$ 0.041  &  0.878  \\
Sin & 0.689  &  0.973  &  0.71  &  0.697 $\pm$ 0.006  &  0.766  \\
Trend & 0.684  &  0.972  &  0.786  &  0.697 $\pm$ 0.006  &  0.743  \\ \hline
\multicolumn{6}{|c|}{100,000 points FARIMA (0,d,0) --- H = 0.9.} \\ \hline
None & 0.757  &  0.882  &  0.91  &  0.886 $\pm$ 0.004  &  0.861  \\
AR(1) & 0.804  &  0.789  &  0.873  &  0.969 $\pm$ 0.036  &  1.011  \\
Sin & 0.764  &  0.967  &  0.913  &  0.886 $\pm$ 0.004  &  0.883  \\
Trend & 0.757  &  0.974  &  0.933  &  0.886 $\pm$ 0.004  &  0.875  \\ \hline
\multicolumn{6}{|c|}{100,000 points FARIMA (1,d,1) --- H= 0.9,
$\phi_1 = 0.5, \theta_1 = 0.5$.} \\ \hline
None & 0.856  &  0.854  &  0.881  &  0.887 $\pm$ 0.006  &  0.858  \\
AR(1) & 0.888  &  0.773  &  0.874  &  0.959 $\pm$ 0.04  &  1.001  \\
Sin & 0.86  &  0.963  &  0.885  &  0.887 $\pm$ 0.006  &  0.879  \\
Trend & 0.856  &  0.968  &  0.92  &  0.887 $\pm$ 0.006  &  0.872  \\ \hline
\multicolumn{6}{|c|}{100,000 points FARIMA (2,d,1) --- H= 0.9, 
$\phi_1 = 0.5, \phi_2 = 0.2, \theta_1 = 0.1$.} \\ \hline
None & 0.807  &  0.74  &  0.817  &  0.966 $\pm$ 0.048  &  1.05  \\
AR(1) & 0.814  &  0.691  &  0.822  &  1.007 $\pm$ 0.059  &  1.136  \\
Sin & 0.8  &  0.94  &  0.821  &  0.966 $\pm$ 0.048  &  1.052  \\
Trend & 0.807  &  0.939  &  0.856  &  0.966 $\pm$ 0.048  &  1.051  \\ \hline
\end{tabular}
\end{center}
\caption{Results for various FARIMA models corrupted by 
several forms of noise.}
\label{tab:farima}
\end{table*}

Table \ref{tab:farima} shows a variety of results for FARIMA models.
The first three runs are for a FARIMA $(0,d,0)$ model (that is one
with no AR or MA components) and with a Hurst parameter $H = 0.7$.
In this case, all methods peform adequately with no noise (although
the R/S plot perhaps underestimates the answer).  Addition of AR(1)
noise causes problems for the R/S plot, wavelet and local Whittle
methods and to a lesser extent the periodogram.  
The addition of a sin wave and a trend causes problems for the aggregated
variance.

For a FARIMA $(1,d,1)$ model with $H = 0.7$ and with the AR parameter
$\phi_1 = 0.5$ and the MA parameter $\theta_1 = 0.5$ (implying a moderate
degree of short range correlation) all estimators provide a reasonable
result for the uncorrupted series.  As before, the wavelet and local
Whittle method seem relatively robust to the addition of a trend.
The AR(1) noise again causes problems for most of the methods.

For a FARIMA $(0,d,0)$ model with $H = 0.9$ the R/S method under
predicts the Hurst parameter but all others perform well in
the absence of noise.  The AR(1) noise causes problems for the
local Whittle and wavelet methods and the sine wave and trend 
cause problems for the aggregated variance.

For a FARIMA $(1,d,1)$ model with $H = 0.9$ and with the AR parameter
$\phi_1 = 0.5$ and the MA parameter $\theta_1 = 0.5$ (implying, as before,
a moderate degree of short range correlation) all estimators do 
relatively well initially.  The corruption produces the same
problems with the same estimators --- that is to say, wavelets and
local Whittle do not cope with the AR(1) noise and Aggregated
variance reacts badly to the sine wave and local trend.

For a FARIMA $(2,d,1)$ model with $H = 0.9$ and with the AR parameters
$\phi_1 = 0.5$, $\phi_2 = 0.2$ and the MA parameter $\theta_1 = 0.1$
indicating quite strong short-range correlations, none of the
estimators perform particularly well.  Presented with these results, 
a researcher would certainly not know
the Hurst parameter of the underlying model from looking at the
results given by the estimators.  In the case of the AR(1) corrupted data the
measurement from the Wavelet estimator is outside of the usual range
for the Hurst parameter.  In fact it is not unusual for Hurst parameter estimators to
produce estimates outside the range $(1/2,1)$.  
All five estimators are producing different
results in most cases (there is some aggreement between the R/S plot and
periodogram but it would be hard to put this down to anything more than 
coincidence and, in any case, they are agreeing on an incorrect
value for the Hurst parameter).  It is interesting that, even in this
relatively simple case where the theoretical correct result is known, five
well-known estimators of the Hurst parameter all fail to get the
correct answer.

\subsection{Autocorrelations for the Artificial Data}

\begin{figure*}
\begin{center}
\includegraphics[width=7.5cm]{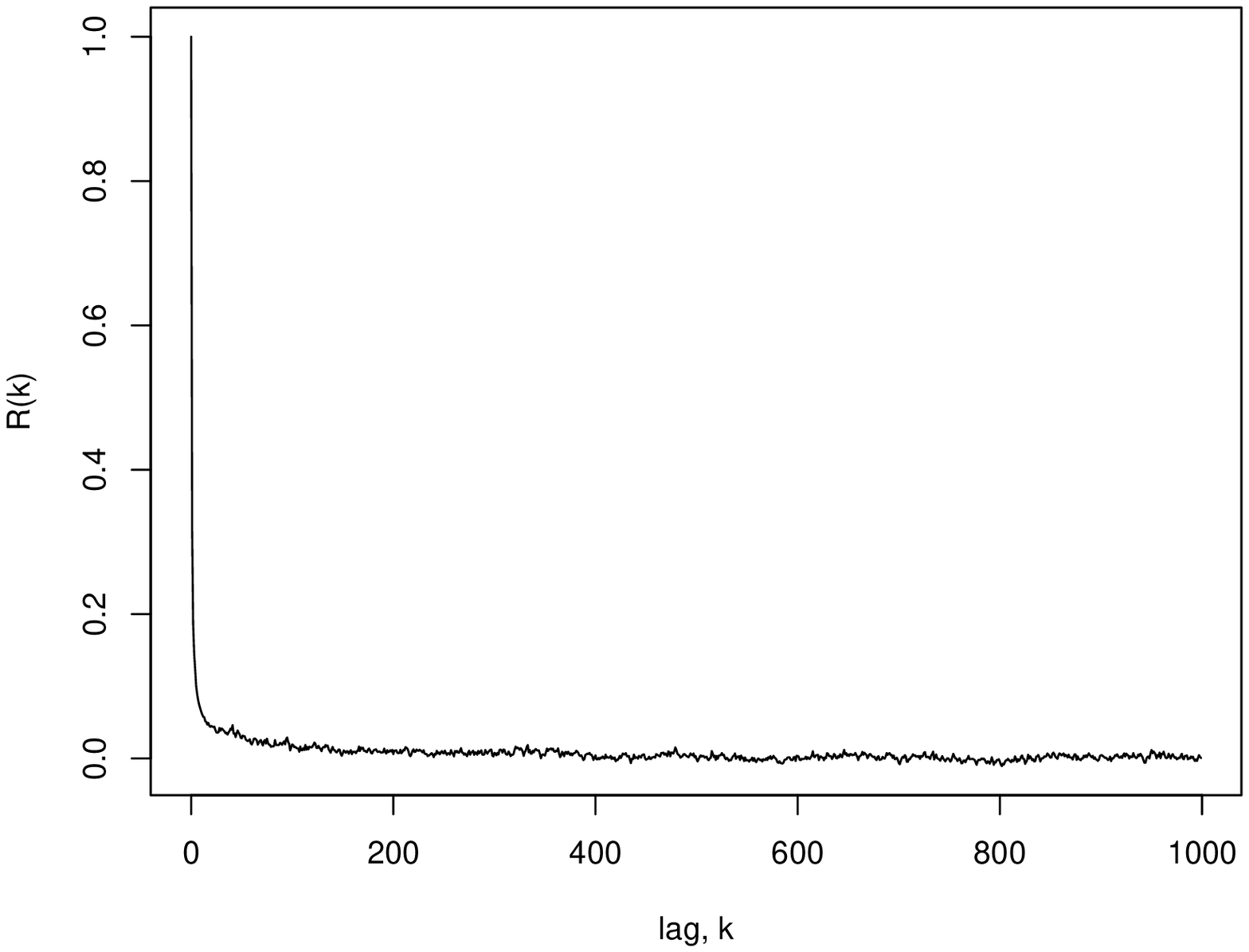}
\includegraphics[width=7.5cm]{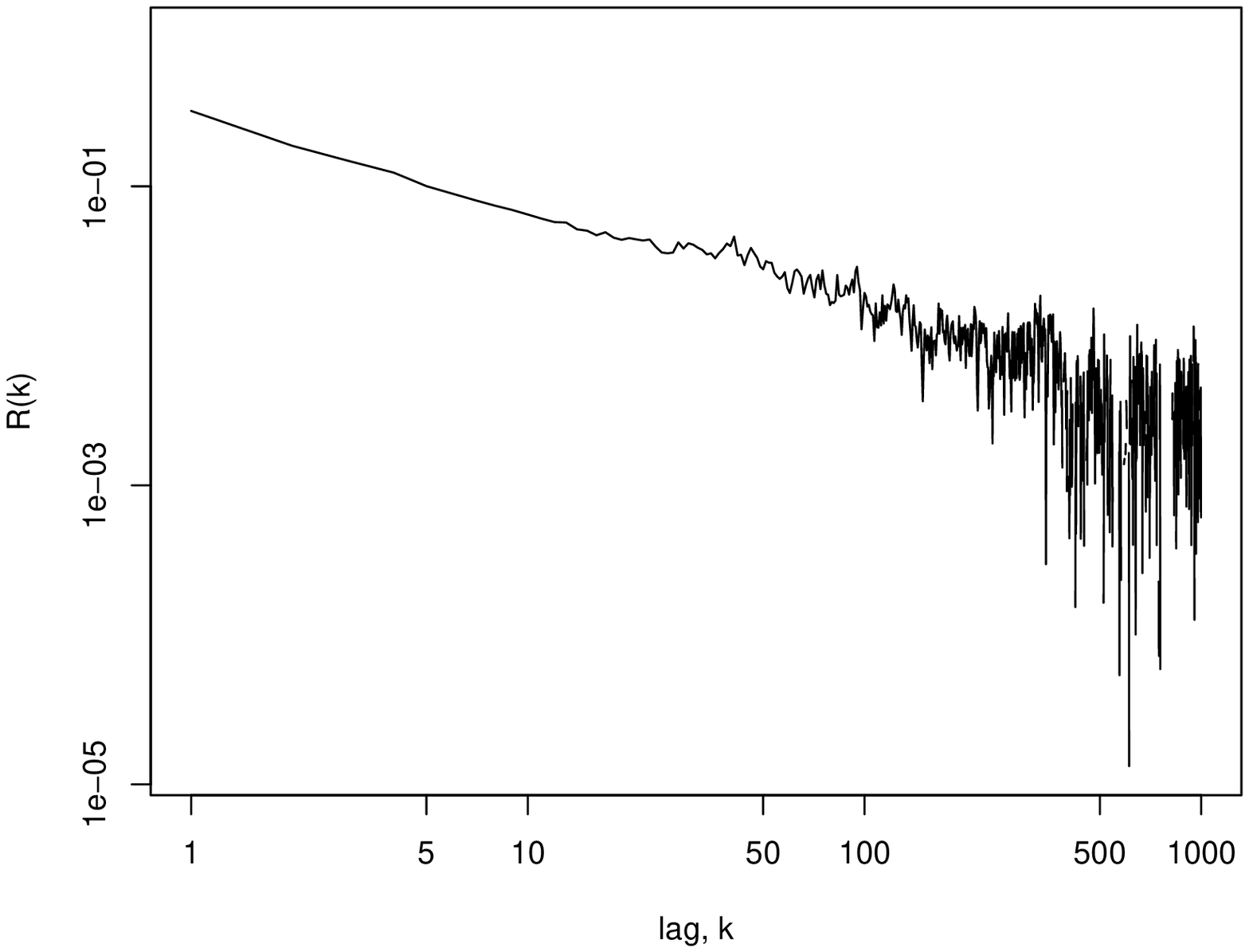}
\end{center}
\caption{ACF (left) and log-log ACF (right) for FGN with Hurst parameter $H=0.7$.}
\label{fig:fgn-acf}
\end{figure*}

It's instructive to look at the ACF of these data sets to understand why
the various methods fail or succeed with the data sets.  Figure 
\ref{fig:fgn-acf} shows the ACF up to lag 1000 for a data set of 100,000 
points of FGN data with $H=0.7$.  For this data, it is possible to fit ``by
eye'' a straight line to the log-log plot of the ACF and obtain an estimate
of the Hurst parameter.  From Table \ref{tab:fgn} it can be seen that all the
estimators performed well on this data set.  Note also that in the log-log plot
it can be seen that at the higher lags the error on the ACF estimate is large.

\begin{figure*}
\begin{center}
\includegraphics[width=7.5cm]{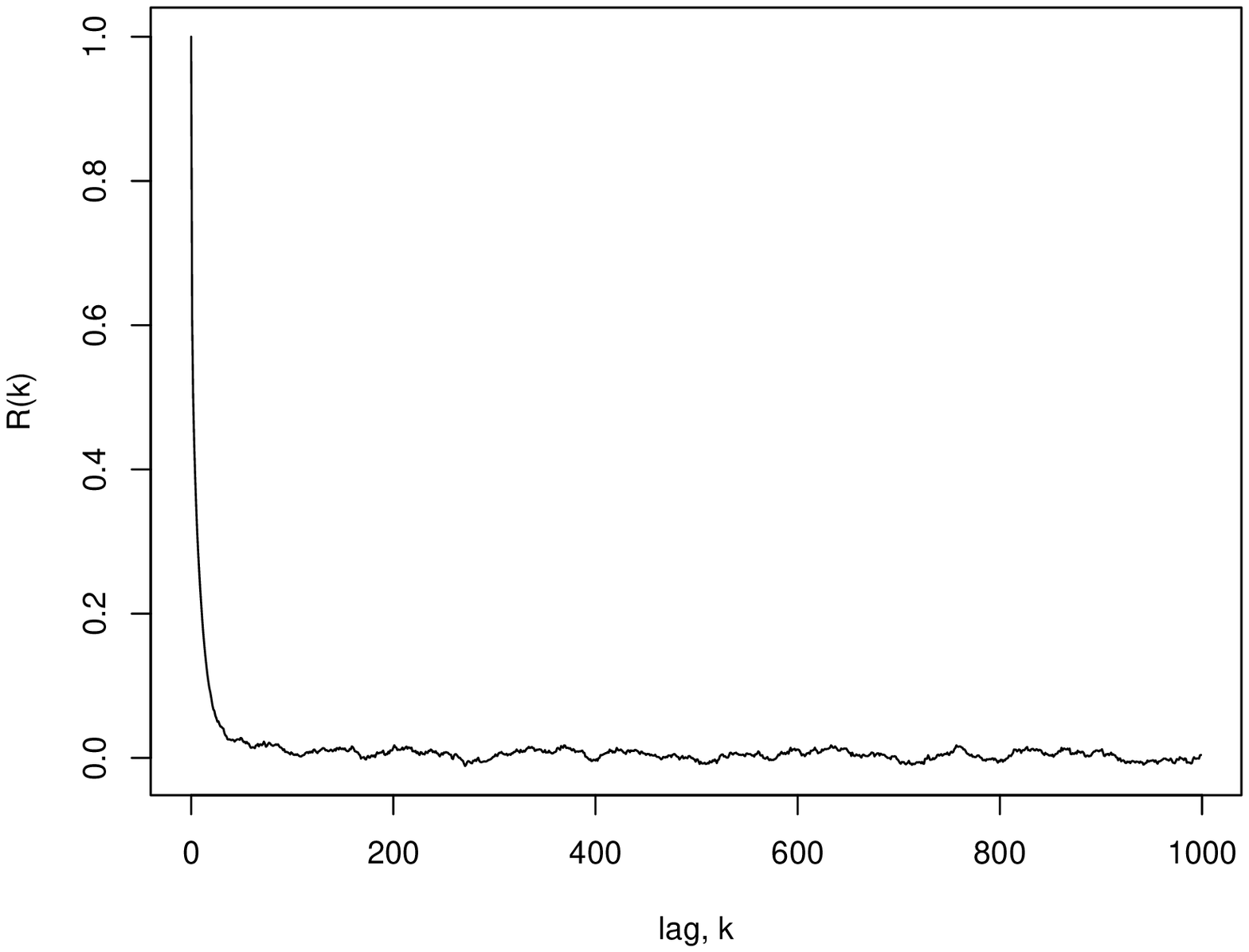}
\includegraphics[width=7.5cm]{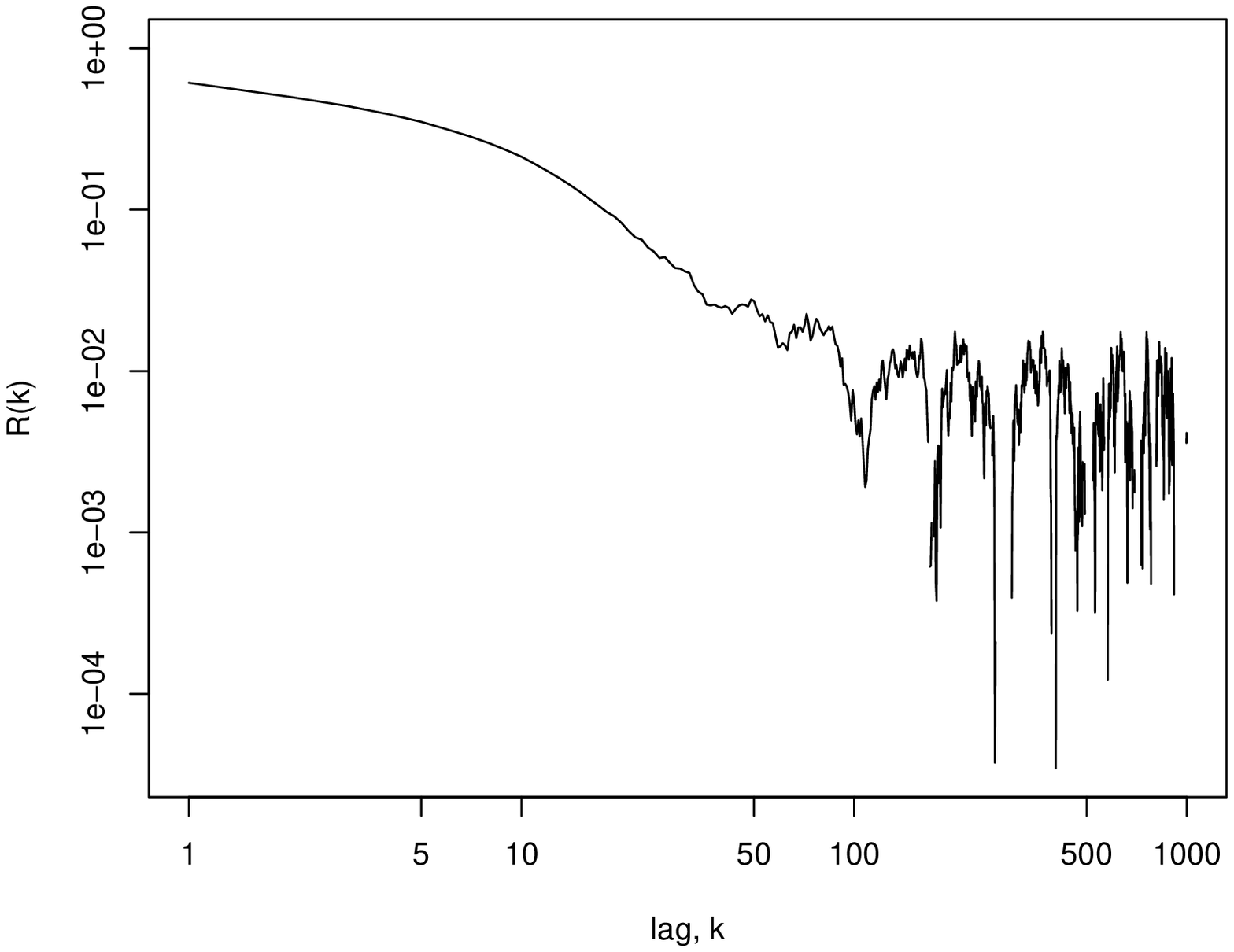}
\end{center}
\caption{ACF (left) and log-log ACF (right) for FGN with Hurst parameter $H=0.7$
corrupted by AR(1) noise with the same variance.}
\label{fig:fgn-acf-noise}
\end{figure*}

When the time series is corrupted by the addition of AR(1) noise as described
earlier in the paper then the ACF changes markedly.  The ACF is then shown
in Figure \ref{fig:fgn-acf-noise}.  The degree to which the ACF has changed 
is only really clear in the log-log plot.  It can be seen that, for low lags,
the ACF remains much higher than in the noise-free data of the previous series.
It would be difficult indeed to make a convincing case for fitting a straight
line to this data.  As for the higher lags, the ACF estimate certainly does not
seem to produce anything like a straight line in the log-log plot for lags over 
fifty.  Two things can be clearly seen from this picture, firstly that it
is impossible to get a good estimate of LRD simply by fitting a straight line
to the ACF and secondly that the addition of highly correlated short range
dependent data can vastly change the nature of the estimation problem.  From
considering this ACF it may be no surprise that the estimators mainly performed
so badly at removing AR(1) noise.

\begin{figure*}
\begin{center}
\includegraphics[width=7.5cm]{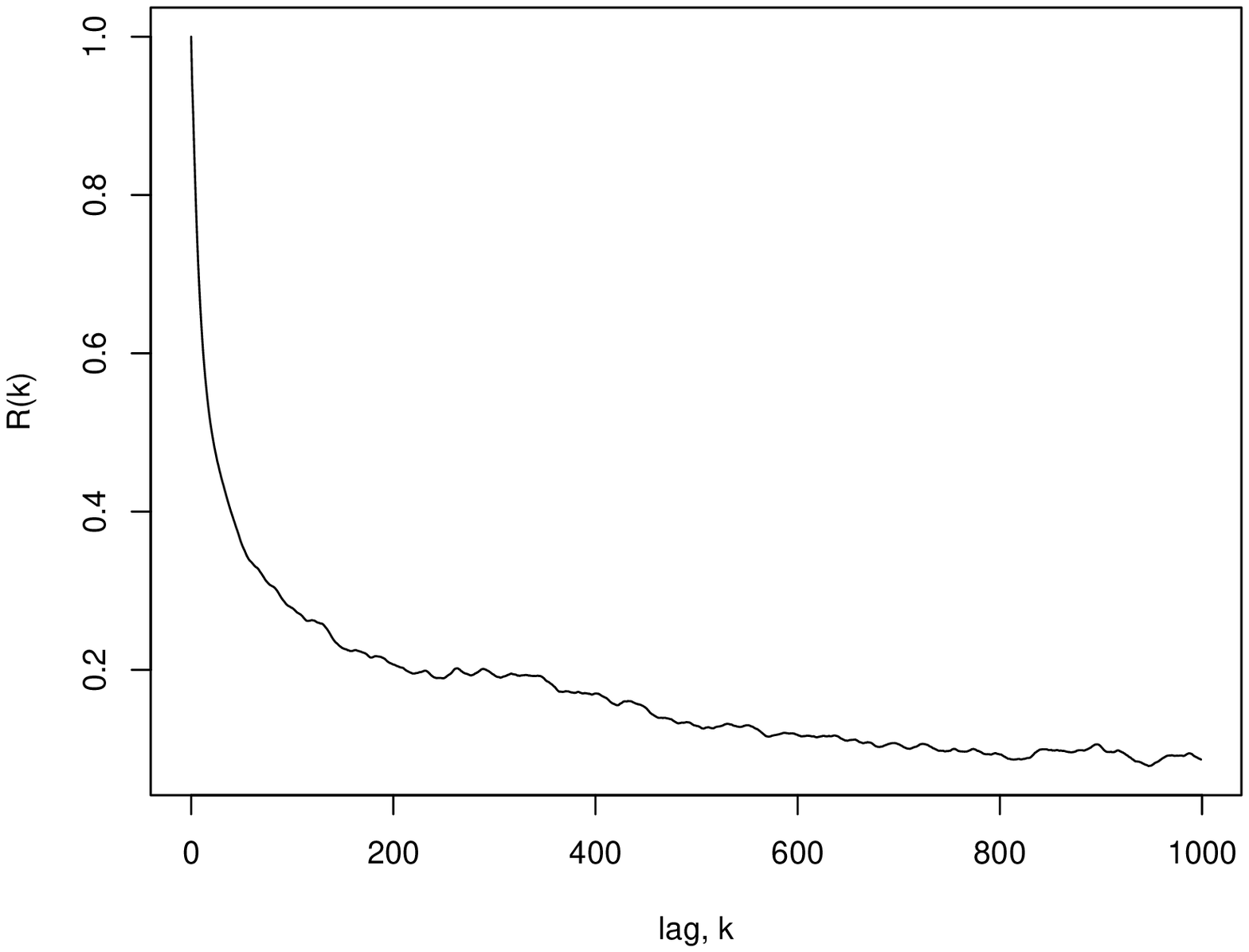}
\includegraphics[width=7.5cm]{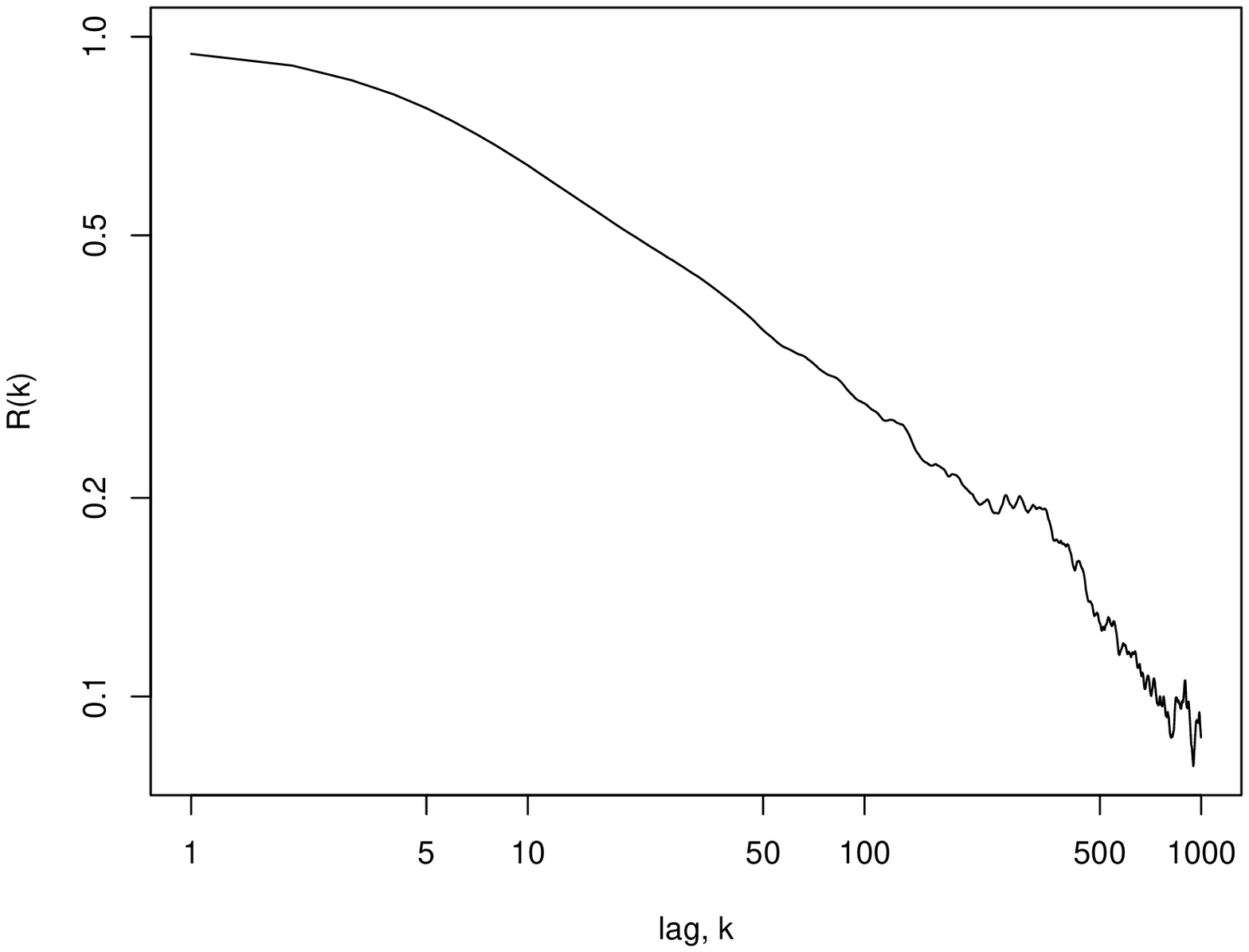}
\end{center}
\caption{ACF (left) and log-log ACF (right) for FARIMA $(2,d,1)$ --- $H= 0.9$, 
$\phi_1 = 0.5, \phi_2 = 0.2, \theta_1 = 0.1$.}
\label{fig:far-acf}
\end{figure*}

Finally, the ACF is shown for the a data set which is FARIMA $(2,d,1)$ with $H= 0.7$, 
$\phi_1 = 0.5, \phi_2 = 0.2, \theta_1 = 0.1$.  This was the data which proved
hardest to estimate in Table \ref{tab:farima}.  Some of the difficulties of
this estimation can be seen by looking at the ACF in Figure \ref{fig:far-acf}.  Even
before the addition of noise, it can be seen that this data looks as hard to find
a single best fit line as it was in Figure \ref{fig:fgn-acf-noise}. It
is, again, unsurprising, that the estimators performed badly with this data
set even without the addition of noise.

\subsection{Results on Real Data}

In analysing the real data it is hard to know where to begin.  Since
the genuine answer (if, indeed, it can be really said that there is a
genuine answer) is not known it cannot be said that one result is more
``right'' than another.  The suggested methods for preprocessing data
(taking logs, removing a linear trend and removing a best fit polynomial
--- in this case of order ten) have all been found in the literature
on measuring the Hurst parameter.

\begin{table*}
\begin{center}
\begin{tabular}{|l| l l l l l|} \hline
 Filter & R/S Plot & Aggreg. & Period. & Wavelet & Local \\
 Type & & Variance & ogram & Estimate & Whittle \\ \hline
\multicolumn{6}{|c|}{York trace (bytes/second) --- 4047 points} \\ \hline
None & 0.749  &  0.88  &  1.186  &  0.912 $\pm$ 0.052  &  0.981  \\
Log & 0.758  &  0.894  &  1.105  &  0.921 $\pm$ 0.039  &  0.932  \\
Trend & 0.749  &  0.873  &  1.212  &  0.912 $\pm$ 0.052  &  0.981  \\
Poly & 0.756  &  0.723  &  0.732  &  0.895 $\pm$ 0.04  &  0.972  \\ \hline
\multicolumn{6}{|c|}{York trace (bytes/tenth) --- 40467 points} \\ \hline
None & 0.826  &  0.924  &  0.928  &  0.909 $\pm$ 0.012  &  0.881  \\
Trend & 0.826  &  0.923  &  0.932  &  0.909 $\pm$ 0.012  &  0.881  \\
Poly & 0.827  &  0.892  &  0.863  &  0.909 $\pm$ 0.012  &  0.878  \\ \hline
\end{tabular}
\caption{Analysis of bytes/unit time data collected at the University of York.}
\label{tab:yorkdata}
\end{center}
\end{table*}

Table \ref{tab:yorkdata} shows analysis of data collected at the
University of York.  The same data set is analysed firstly as a
series of bytes/second and then as bytes/tenth of a second.  While
theoretically the results should be the same, in practice this is
not the case.  Obviously there are only one tenth as many points
in the data set when seconds are used rather than tenths of
seconds.  Firstly, looking at the data aggregated
over a time period of one second, there is no good agreement between
estimators.  The periodogram estimate is hopelessly out of the correct
range.  The other estimators, while in the range $(1/2,1)$ show no 
particular agreement.  Of the suggested filtering techniques, little changes
between them except that removal of a polynomial greatly reduces the
estimate found by the periodogram and slightly reduces the estimate
found by aggregated variance.  No conclusion can realistically be
drawn about the data from these results.

Considering the data aggregated into tenths of a second time units the
picture is somewhat clearer.  Taking a log of data was impossible at
this time scale due to presence of zeros.  The estimators, with the
exception of the R/S plot are all relatively near $H = 0.9$.  While
it seems somewhat arbitrary to ignore the results of the R/S plot it
should be remembered that this technique performed poorly with high
Hurst parameter measurements on theoretical data and underestimated
badly in those cases.  No great difference is observed from any
of the suggested filtering techniques except, perhaps, a slight
reduction in the aggregated variance and periodogram results from
removal of a polynomial.  A tentative conclusion from this data
would be that $0.85 < H < 0.95$ and that the R/S plot is inaccurate
for this trace.

\begin{table*}
\begin{center}
\begin{tabular}{|l| l l l l l|} \hline
 Filter & R/S Plot & Aggreg. & Period. & Wavelet & Local \\
 Type & & Variance & ogram & Estimate & Whittle \\ \hline
\multicolumn{6}
{|c|}{Bellcore data BC-Aug89 (interarrival times) --- first 360,000 points.} \\ \hline
None & 0.73  &  0.742  &  0.762  &  0.73 $\pm$ 0.018  &  0.661  \\
Log & 0.722  &  0.806  &  0.797  &  0.77 $\pm$ 0.02  &  0.652  \\
Trend & 0.73  &  0.74  &  0.762  &  0.73 $\pm$ 0.018  &  0.661  \\
Poly & 0.73  &  0.733  &  0.751  &  0.73 $\pm$ 0.018  &  0.66  \\ \hline
\multicolumn{6}
{|c|}{Bellcore data BC-Aug89 (interarrival times) --- second 360,000 points.} \\ \hline
None & 0.709  &  0.703  &  0.742  &  0.746 $\pm$ 0.025  &  0.655  \\
Log & 0.721  &  0.795  &  0.779  &  0.778 $\pm$ 0.011  &  0.673  \\
Trend & 0.709  &  0.703  &  0.742  &  0.746 $\pm$ 0.025  &  0.655  \\
Poly & 0.709  &  0.691  &  0.732  &  0.746 $\pm$ 0.025  &  0.654  \\ \hline
\multicolumn{6}
{|c|}{Bellcore data BC-Aug89 (bytes/10ms) --- first 1000 secs.} \\ \hline
None & 0.707  &  0.8  &  0.817  &  0.786 $\pm$ 0.017  &  0.822  \\
Trend & 0.707  &  0.797  &  0.815  &  0.786 $\pm$ 0.017  &  0.822  \\
Poly & 0.707  &  0.789  &  0.787  &  0.786 $\pm$ 0.017  &  0.822  \\ \hline
\multicolumn{6}
{|c|}{Bellcore data BC-Aug89 (bytes/10ms) --- second 1000 secs.} \\ \hline
None & 0.62  &  0.802  &  0.808  &  0.762 $\pm$ 0.012  &  0.825  \\
Trend & 0.62  &  0.802  &  0.808  &  0.762 $\pm$ 0.012  &  0.825  \\
Poly & 0.618  &  0.786  &  0.777  &  0.762 $\pm$ 0.012  &  0.824  \\ \hline
\end{tabular}
\caption{Analysis of bytes/unit time and interarrival times for the Bellcore
data with various methods to attempt to remove non-stationary components.}
\end{center}
\end{table*}

In the case of the Bellcore measurements, the data has been split into
two sections and analysed seperately for interarrival times and 
for bytes per unit time.  Considering first the interarrival times,
all estimators seem to have a result which is not too distant from
$H = 0.7$ in both cases.  The various filtering techniques tried
do little to change this.  It is hard to come to a really robust
conclusion since the estimators are as high as $0.806$ (aggregated
variance after taking logs) and as low as $0.652$ (local Whittle
after taking logs).

When the bytes per unit time are considered, the log technique cannot
be used due to zeros in the data.   The most comfortable conclusion
abou this data might be that the Hurst parameter is somewhere around
$H = 0.8$ with the R/S plot underestimating again.  As before, it
is hard to reach a strong conclusion on the exact Hurst parameter.
Certainly it would be foolish to take the confidence intervals
on the wavelet estimator at face value.  The various filters tried
seem to make little difference except perhaps a slight reduction 
in the answer given by some estimators after the polynomial is removed.
A tentative conclusion might be that $0.75 < H < 0.85$ for this
data with the R/S plot being in error.

\section{CONCLUSION}

This paper has looked at measuring the Hurst parameter, firstly in
the case of artificial data contaminated by various types of noise and
secondly in the case of real data with various filters to try to
improve the performance of the estimators used.  

The most striking conclusion of this paper is that measuring the Hurst
parameter, even in artificial data, is very hit and miss.  In
the artificial data with no corrupting noise, some estimators performed
very poorly indeed.  Confidence intervals given should certainly not
be taken at face value (indeed should be considered as next to worthless).

Corrupting noise can affect the measurements badly and different estimators
are affected in by different types of noise.  In particular, frequency
domain estimators (as might be expected) are robust to the addition
of sinusoidal noise or a trend.  All estimators had problems in some
circumstances with the addition of a heavy degree of short-range dependence
even though this, in theory, does not change the long-range dependence
of the time series.

When considering real data, researchers are advised to use extreme caution.
A researcher relying on the results of any single estimator for the 
Hurst parameter is likely to be drawing false conclusions, no matter how
sound the theoretical backing for the estimator in question.  While simple 
filtering techniques are suggested in the literature for improving the
performance of Hurst parameter estimation, they had little or no
effect on the data analysed in this paper.

All the data and tools used in this paper are available for download
from the web and can be found at: \\
{\tt www.richardclegg.org/ \\ lrdsources/software/}
\bibliography{rgc_ijs}

\section*{BIOGRAPHY}

\textbf{ Richard G. Clegg} is a Research Assistant in the Department of Mathematics at the
University of York.  He obtained his PhD in ``The Statistics of Dynamic Networks'' in
May 1994 and works on applied mathematics in networks, mainly road networks and
traffic networks.  His main research interests are equilibrium and driver route
choice in road networks and long-range dependence in computer networks.

\end{document}